\documentclass[12pt]{amsart}
\usepackage{amsfonts,amssymb,enumerate,color,url,graphicx, epsfig}

\input xy
\xyoption{all}

\def\Pr{\begin{proof}}
\def\Rp{\end{proof}}

\def\norm{\left\Vert{.}\right\Vert}

\newcommand\norme[1]{\|#1\|}
\def\IR{\mathbb R}
\def\IN{\mathbb N}

\def\ZF{\mathbf {ZF}}
\def\ZFC{\mathbf {ZFC}}

\def\DC{\mathbf {DC}}
\def\HB{\mathbf {HB}}

\def\BPI{\mathbf {BPI}}

\def\AC{\mathbf {AC}}

\def\ACD{\mathbf {AC(\IN)}}

\def\ACDF{\mathbf {AC(\IN,fin)}}

\def\ZFC{\mathbf {ZFC}}

\def\AUc{\mathbf {A1}}
\def\AH{\mathbf {A2}}
\def\AHb{\mathbf {A3}}
\def\AHbf{\mathbf {A4}}

\DeclareMathOperator{\spanv}{span}
\DeclareMathOperator{\conv}{conv}
\DeclareMathOperator{\diam}{diam}

\theoremstyle{plain}

\newtheorem{corollary}{Corollary}
\newtheorem{proposition}{Proposition}
\newtheorem*{proposition*}{Proposition}
\newtheorem{theorem}{Theorem}
\newtheorem*{theorem*}{Theorem}
\newtheorem*{corollary*}{Corollary}
\newtheorem{lemma}{Lemma}
\newtheorem*{lemma*}{Lemma}

\theoremstyle{definition}
\newtheorem{definition}{Definition}
\newtheorem{notation}{Notation}
\newtheorem{question}{Question}

\theoremstyle{remark}
\newtheorem{remark}{Remark}
\newtheorem{example}{Example}

\date{\today}

\begin{document}
\title{Countable choice and compactness}
\author[M.~Morillon]{Marianne Morillon}
 \address{ERMIT, D\'epartement de Math\'ematiques et Informatique,
 Universit\'e de La R\'eunion, 15 avenue Ren\'e Cassin - BP 7151 -
 97715 Saint-Denis Messag. Cedex 9 FRANCE}
 \email[Marianne Morillon]{mar@univ-reunion.fr}
 \urladdr{http://personnel.univ-reunion.fr/mar}

\begin{abstract} 
We work in set-theory without choice $\ZF$.
Denoting by  $\ACD$ the countable axiom of choice, we show in $\ZF \mathbf + \ACD$ 
that the closed unit ball of a uniformly convex Banach space is compact in the convex 
topology (an alternative to the weak topology in $\ZF$). We prove that 
this  ball is (closely) convex-compact in the convex topology. 
Given a set $I$, a real number $p \ge 1$ ({\em resp.}  $p=0$), and some closed subset $F$ 
of $[0,1]^I$ which is a bounded subset
of $\ell^p(I)$, we show that $\ACD$ ({\em resp.} $\DC$, the axiom of Dependent Choices) 
implies the compactness of $F$.
\end{abstract}

 \maketitle

 \begin{center}
    \mbox{\epsfbox{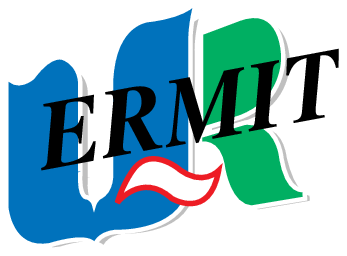}}
    \end{center}
 \begin{center}
 EQUIPE R\'EUNIONNAISE DE MATH\'EMATIQUES ET INFORMATIQUE TH\'EORIQUE (ERMIT)
 \end{center}
 \bigskip
%
%
%
 \begin{center}
 http://laboratoires.univ-reunion.fr/ermit
 \end{center}
%
%
%
%
%

\tableofcontents

\section{Introduction}
\subsection{Presentation of the results}
We work in  $\ZF$, Zermelo-Fraenkel set-theory without the Axiom of Choice (for short $\AC$). Consider  the {\em countable Axiom of Choice}, which is not provable in $\ZF$, and which does not imply $\AC$:
\begin{quote}$\ACD$: {\em If $(A_n)_{n \in \IN}$ is a family of non-empty sets, then there exists a mapping $f : \IN \to \cup_{n \in \IN}A_n$ associating to every
$k \in \IN$ an element $f(k) \in A_k$.}
\end{quote}

In this paper, we  first provide in $\ZF \mathbf + \ACD$ a criterion of compactness 
for  topological spaces   coarser  than  some   complete metric space, having a sub-basis of closed sets   satisfying ``good'' properties with respect to  
the distance   (see Theorem~\ref{theo:small-subbasis2compact} in Section~\ref{subsubsec:criterion}).  We then  consider an alternative topology for the weak topology on a normed space, namely  the {\em convex topology} (which  is the weak topology in $\ZF \mathbf + \HB$), and we provide  some properties of this convex topology: in particular, using the (choiceless) Lusternik-Schnirelmann theorem, we show (see Theorem~\ref{theo:clos-ball} in Section~\ref{subsec:convex-top}) that the closure of the unit sphere of $E$ for the convex topology is the closed unit ball of $E$. Applying  our  criterion of compactness to the convex topology, we obtain some  new  results.
First, in $\ZF \mathbf + \ACD$,  ``The closed unit ball of a uniformly convex Banach space is compact in the convex topology.'' 
(see Theorem~\ref{theo:acd2wc-uc} in Section~\ref{subsec:weak-alaoglu}): this extends a result obtained by  Fremlin for Hilbert spaces (see \cite[Chapter~56, Section 566P]{Fr06}) and this solves a question raised
in \cite[Question~2]{De-Mo}. We then prove in $\ZF$  that ``The closed unit ball of a uniformly convex Banach space is 
(closely) {\em convex-compact} in the convex topology.'' (see Theorem~\ref{theo:uc-ccomp-ZF} in Section~\ref{subsec:cc-dual-uc}).  Given a set $I$, we apply our results to closed subsets of $[0,1]^I$. 
In Section~\ref{sec:acd-comp}, we show that $\ACD$ implies the compactness of 
closed subsets of $[0,1]^I$ which are bounded subsets
of some $\ell^p(I)$, $1 \le p <+\infty$. In Section~\ref{sec:dc-and-comp}, we prove that the Axiom of {\em Dependent Choices}
$\DC$ implies that every closed subset of $[0,1]^I$ which is contained in $\ell^0(I)$ is compact. 

\subsection{Some weak forms of $\AC$}
We now review some weak forms of the Axiom of Choice  which  will be used in this paper and the known links between them. 
For detailed references and much  information on this subject, see \cite{Ho-Ru}.

\subsubsection{$\DC$ and $\ACDF$} The axiom of {\em Dependent Choices} asserts that:
\begin{quote} $\DC$: Given a non-empty set $X$ and a binary relation $R$ on $X$ such that $\forall x \in X \exists y \in X \; xRy$, there exists a sequence 
$(x_n)_{n \in \IN}$ of $X$ such that for every $n \in \IN$, $x_n R x_{n+1}$.
\end{quote} 
The {\em countable Axiom of Choice for finite sets} says that:
\begin{quote}$\ACDF$: {\em If $(A_n)_{n \in \IN}$ is a family of finite non-empty sets, then there exists a mapping $f : \IN \to \cup_{n \in \IN}A_n$ associating to every
$n \in \IN$ an element $f(n) \in A_n$.}
\end{quote} 

Of course, $\AC \Rightarrow \DC \Rightarrow \ACD \Rightarrow \ACDF$. However, the converse statements are not provable in $\ZF$, and $\ACDF$ is not provable in $\ZF$ (see references in \cite{Ho-Ru}).


\subsubsection{$\BPI$ and $\HB$} The {\em Boolean Prime Ideal} axiom says that:
\begin{quote} $\BPI$: {\em Every non trivial boolean algebra has a prime ideal.}
\end{quote}
It is known that $\BPI$ is not provable in $\ZF$ and that $\BPI$ does not imply $\AC$. The following well-known statements of functional analysis 
are equivalent to the axiom $\BPI$ (see for example \cite{Sa-Bo}): the {\em Tychonov theorem} for product of compact Hausdorff spaces, 
the {\em Alaoglu theorem}, the fact that for every set $I$ the product space $[0,1]^I$ is compact.

\begin{remark} \label{rem:tycho-wo} If a set $I$ is well-orderable, then the product topological space $[0,1]^I$ is compact in $\ZF$.
\end{remark}

\subsubsection{Hahn-Banach} \label{subsubsec:HB}
Given a (real) vector space $E$, a  mapping $p: E \to \IR$ is {\em sub-linear} if 
 for every $x, y \in E$, and every $\lambda \in \IR_+ $, $p(x+y) \le p(x) + p(y)$ (sub-additivity),
 and $p(\lambda.x) = \lambda p(x)$ (positive homogeneity). Consider the  {\em ``Hahn-Banach''}  axiom, a well known consequence  of $\AC$ which is not provable in $\ZF$:
\begin{quote}$\HB$: {\em Let $E$ be a (real) vector space.  
If  $p:E \to \IR$ is a sub-linear mapping, if $F$ is a vector subspace of $E$, if  $f: F \to \IR$ is a linear mapping such that $f \le p_{\restriction F}$, then there exists 
a linear mapping $g : E \to \IR$ extending $f$  such that $g \le p$. }
\end{quote}
Given a (real) topological vector space $E$ ({\em i.e.} $E$ is a vector space such that the 
``sum'' $+ : E \times E \to E$ and the external multiplicative law $. : \IR \times E \to E$ are continuous for the product topology),
 say that $E$ satisfies  the  {\em Continuous Hahn-Banach property} (for short {\em CHB} property) if  
{\em ``For every {\em continuous} sub-linear mapping
  $p:E \to \IR$, for every vector subspace $F$ of $E$, if  $f: F \to \IR$ is a linear mapping such that 
$f \le p_{\restriction F}$, then there exists 
a linear mapping $g : E \to \IR$ extending  $f$  such that $g \le p$.''} 
The statement $\HB$  is not provable in $\ZF$, however, some real normed spaces  satisfy (in $\ZF$) the  {\em CHB} property: 
for example normed spaces with a well-orderable dense subset (in particular  separable normed spaces), but also
Hilbert spaces, spaces $\ell^0(I)$ 
(see \cite{Fo-Mo}),  uniformly convex Banach spaces with a G\^ateaux\--differentiable norm (\cite{Do-Mo}), 
uniformly smooth Banach spaces (see~\cite{Al-Mo}).  

It is rather easy to prove that $\BPI$ implies  $\HB$ 
and that $\HB$   is equivalent to most of its classical geometrical forms (see \cite[Section~6]{Do-Mo}).
  It is also easy to see that $\AC \Rightarrow (\BPI \mathbf + \DC) \Rightarrow \BPI \Rightarrow \HB$. 
The converse statements are not
provable in $\ZF$ and $\HB$ is not provable in $\ZF \mathbf + \DC$ (see \cite{Ho-Ru}). 

\begin{remark} There exist various ($\ZFC$-equivalent) definitions of reflexivity for Banach spaces: most of them are equivalent in
$\ZF \mathbf + \HB \mathbf + \DC$ (see \cite{Mo04}), including  James' sup theorem  (see \cite{Mo05}).
\end{remark}

\section{A criterion of compactness}
\subsection{Filters}
\subsubsection{Filters in lattices of sets}
Given a set $X$, a {\em lattice} of subsets of $X$ is a subset $\mathcal L$ of $\mathcal P(X)$ containing $\varnothing$
and $X$, which is closed by finite intersections and finite unions. A {\em filter} of the lattice $\mathcal L$ is a 
non-empty proper subset $\mathcal F$ of $\mathcal L$ such that for every $A, B \in \mathcal L$:
\begin{equation} 
(A, B \in \mathcal F) \Rightarrow A \cap B \in \mathcal F
\end{equation}
\begin{equation} 
(A \in \mathcal F \text{ and } A \subseteq B) \Rightarrow B \in \mathcal F
\end{equation}
A subset $\mathcal A$ of $\mathcal L$ is contained in a filter of $\mathcal L$ if and only if $\mathcal A$ is 
{\em centered} {\em i.e.}   every finite subset 
of $\mathcal A$ has a non-empty intersection;  in this case,
the intersection of all filters of $\mathcal L$ containing $\mathcal A$ is called the {\em filter generated by $\mathcal A$}
and we denote it by $fil(\mathcal A)$.

\subsubsection{Stationary sets} \label{subsubsec:stat}  
Given a filter $\mathcal F$ of a lattice $\mathcal L$ of subsets of a set $X$, an element $S \in \mathcal L$ is 
{\em $\mathcal F$-stationary} if for every $A \in \mathcal F$, $A \cap S \neq \varnothing$. 
The set $\mathcal S_{\mathcal L}(\mathcal F)$ (also denoted by $\mathcal S(\mathcal F)$) of $\mathcal F$-stationary elements of $\mathcal L$ satisfies the following properties:
\begin{enumerate}[(i)]
\item  \label{it:stat1} If $\mathcal A$ is a chain of $\mathcal L$ and if $\mathcal A \subseteq \mathcal S(\mathcal F)$, then 
$\mathcal A \cup \mathcal F$ is centered. 
\item   \label{it:stat2} 
Let $F_1, \dots, F_m \in \mathcal L$.
If $F_1 \cup \dots \cup F_m \in \mathcal S(\mathcal F)$, then there exists some $i_0 \in \{1..m\}$
such that   $F_{i_0}$ is  $\mathcal F$-stationary. 
\end{enumerate}

 \subsection{Complete metric spaces}
Given a metric space  $(X,d)$, some point $a \in X$ and  real numbers $R,R'$ satisfying  $R \le R'$, 
 we define  
 {\em large $d$-balls} and {\em large $d$-crowns} as follows: 
$$B_d(a,R):= \{x \in X : \; d(a,x)\le R\}$$ 
$$D_d(a,R,R') := \{x \in X : \; R \le d(a,x) \le R'\}$$ 
Moreover, if $A$ is a  
subset of $X$, we define the {\em $d$-diameter} of $A$:
$$\diam_{d}(A) := \sup \{d(x,y) : x, y \in A\} \in [0,+\infty]$$
In particular, $\diam_d(\varnothing)=0$.
 A metric space $(X,d)$ is said to be {\em complete} if  every Cauchy filter of the lattice of closed subsets of $X$ has a non-empty intersection. 
Here, a set $\mathcal A$ of subsets of $X$ is {\em Cauchy} if for every $\varepsilon >0$, there exists $A \in \mathcal A$
satisfying $\diam_d(A) < \varepsilon$.

\subsection{A criterion of compactness in  $\ZF \mathbf + \ACD$}
\subsubsection{Compactness}
\begin{definition}[$\mathcal C$-compactness, closed $\mathcal C$-compactness] 
Given a class  $\mathcal C$ of  subsets of a set $X$, say that a subset 
$A$ of $X$ is {\em $\mathcal C$-compact} if for every family $(C_i)_{i \in I}$ of $\mathcal C$ such that
$(C_i \cap A)_{i \in I}$ is centered, $A \cap \cap_{i \in I}C_i$ is non-empty;
say that $\mathcal A$ is {\em closely} $\mathcal C$-compact if there is a mapping associating 
to  every family $(C_i)_{i \in I}$ of $\mathcal C$ such that
$(C_i \cap A)_{i \in I}$ is centered, an element of $A \cap \cap_{i \in I}C_i$.
\end{definition}

Recall that a topological space $X$ is {\em compact} if $X$ is $\mathcal C$-compact, where $\mathcal C$ is the set of closed subsets of $X$.
 Equivalently, every filter of the lattice of closed sets of $X$ has
a non-empty intersection. 

\subsubsection{Sub-basis of closed sets}
\begin{definition}[basis, sub-basis of closed subsets]
A set $\mathcal B$ of closed subsets of a topological space $X$  is a {\em basis of closed sets} if
every closed set of $X$  is an intersection of elements of $\mathcal B$.
A set $\mathcal S$ of closed subsets of $X$   is a {\em sub-basis of closed sets} if
the set $\mathcal B$ of finite unions of elements of $\mathcal S$ is a basis of closed sets of $X$.
\end{definition}

The following result is easy.
\begin{proposition} \label{prop:crit-comp} Let $X$ be a topological space,  
and $\mathcal L$ be a lattice of closed subsets of $X$.  If $\mathcal L$ is a basis of closed subsets of $X$,
and if every filter of $\mathcal L$ 
has a non-empty intersection, then $X$ is compact.   
\end{proposition}

\subsubsection{Property of smallness}
Given real numbers $a,b$, we denote by $]a,b[$ the open interval $\{x \in \IR: \; a<x<b\}$. 
\begin{definition}[smallness in thin crowns] Let $(X,d)$ be a metric space and let $a \in X$.  
Say that a set $\mathcal C$ of subsets of $X$ satisfies the {\em property of $d$-smallness 
in thin crowns centered at $a$} if for every $R \in \IR^*_+$, for every $\varepsilon >0$ there exists $\eta \in ]0,R[$ such that for every $C \in \mathcal C$,
$$C \subseteq D_{d}(a,R-\eta,R+\eta) \Rightarrow \diam_{d}(C)<\varepsilon$$
\end{definition}

\subsubsection{Criterion of compactness} \label{subsubsec:criterion}
\begin{theorem} \label{theo:small-subbasis2compact} Let $(X,d)$ be a complete metric space. 
Let $\mathcal T$ be a topology on $X$ which is included in the topology $\mathcal T_d$ of $(X,d)$.
Let $\mathcal C$ be    a sub-basis  of closed sets of $(X,\mathcal T)$, which is closed by finite intersection.
If  $a \in X$,  if  $\mathcal C$  contains all large $d$-balls centered at $a$, and if 
 $\mathcal C$ satisfies the property of $d$-smallness in thin crowns centered at $a$, then:
\begin{enumerate}[(i)]
\item \label{it:main-ACD} In $\ZF \mathbf + \ACD$,  every large $d$-ball with center $a$ is 
$\mathcal T$-compact  (and thus, every  $d$-bounded 
$\mathcal T$-closed subset of $X$ is $\mathcal T$-compact).
\item \label{it:main-ZF} In $\ZF$, every  large $d$-ball with center $a$ is $\mathcal C$-compact (and thus, every $d$-bounded element of $\mathcal C$  is  closely $\mathcal C$-compact). 
\end{enumerate}
\end{theorem}
\Pr  Let $\mathcal L$ be the sub-lattice of $\mathcal P(X)$ generated by $\mathcal C$. 
Let  $\rho>0$ and let $B$ be the large $d$-ball $B_d(a,\rho)$. \\
\eqref{it:main-ACD}   Let  $\mathcal F$ be a filter of  $\mathcal L$ containing $B$. 
Let us prove in $\ZF \mathbf + \ACD$ that  $\cap \mathcal F$ is non-empty (using  Proposition~\ref{prop:crit-comp}, this will imply that 
$B$ is $\mathcal T$-compact). Let $R:= \inf \{r \in \IR_+ : \; B_d(a,r) \in \mathcal S(\mathcal F)\}$.  
The set of balls $\{B_d(a,r) : \; r >R\}$ is a chain of $\mathcal F$-stationary sets of $\mathcal L$,  thus  
$\mathcal F \cup \{B_d(a,r) : \; r >R\}$ generates a filter  $\mathcal G$ of $\mathcal L$ (see Section~\ref{subsubsec:stat}-\eqref{it:stat1}).  
If  $R=0$ then  $\cap \mathcal G=\{a\}$ (because elements of $\mathcal F$ are $\mathcal T_d$-closed) thus 
$a \in \cap \mathcal G \subseteq \cap \mathcal F$.  
If $R>0$, for every $\varepsilon >0$, there exists some element of  $ \mathcal G$ which is included in the  
crown $D_d(a,R-\varepsilon,R+\varepsilon)$; with $\ACD$, choose for every  $n \in \IN$, a finite subset  $\mathcal Z_n$ of   $\mathcal C$
such that $\cup \mathcal Z_n \in \mathcal G$ and  $\cup \mathcal Z_n \subseteq D_d(a,R-\frac{1}{n+1},R+\frac{1}{n+1})$.  
With $\ACDF$, the set $\cup_{n \in \IN} \mathcal Z_n$ is countable.  We define by induction a sequence 
$(C_n)_{n \in \IN} \in \prod_{n \in \IN} \mathcal Z_n$ such that for every  $n \in \IN$, 
$\mathcal G \cup \{C_i : \; i<n\}$ generates a filter  $\mathcal G_n$ and  
$C_n \in \mathcal S(\mathcal G_n)$: given some $n \in \IN$, 
$\cup \mathcal Z_n \in \mathcal G \subseteq fil(\mathcal G, (C_i)_{i<n}) \subseteq \mathcal S(fil(\mathcal G, (C_i)_{i<n}))$; using Section~\ref{subsubsec:stat}-\eqref{it:stat2}, it follows that there exists 
$C_n \in \mathcal Z_n$ satisfying $C_n \in \mathcal S(fil(\mathcal G, (C_i)_{i<n}))$.   Since $\mathcal C$ satisfies the property of  $d$-smallness in thin $d$-crowns centered at $a$,   the filter  $\mathcal H:=\cup_{n \in \IN} \mathcal G_n$ is  Cauchy in the metric  space 
$(X,d)$.   Since $\mathcal T \subseteq \mathcal T_d$ and $(X,d)$ is complete, 
$\cap \mathcal H$ is a singleton $\{b\}$. 
Thus $b \in   \cap \mathcal H \subseteq \cap \mathcal G \subseteq \cap \mathcal F$.\\
\eqref{it:main-ZF}  Let $\mathcal A$ be  subset of non-empty elements $\mathcal C$ which is closed by finite intersection and such that the ball $B$ belongs to $\mathcal A$. Let us show (in $\ZF$) that $\cap \mathcal A$ is non-empty. 
Let $\mathcal F$ be the filter of $\mathcal L$ generated by $\mathcal A$. Let $R:= \inf \{r \in \IR_+ : \; B_d(a,r) \in \mathcal S(\mathcal F)\}$.  
Denote by $\mathcal A'$ the set $\{A \cap B_d(a,r) : A \in  \mathcal A \text{ and }  r >R\}$.
If  $R=0$ then  $\{a\} = \cap \mathcal A' \subseteq  \cap \mathcal A$.  If $R>0$, then 
for every $\varepsilon >0$, there exists some element of  $ \mathcal A'$ which is included in the  
crown $D_d(a,R-\varepsilon,R+\varepsilon)$.    Since $\mathcal C$ satisfies the property of $d$-smallness 
in thin $d$-crowns centered at $a$,
the centered family  $\mathcal A'$ is  Cauchy  in the metric  space $(X,d)$; since this metric  space  $(X,d)$ is complete,  
 $\cap \mathcal A'$ is a singleton $\{b\}$,  and $\{b\} = \cap \mathcal A' \subseteq \cap \mathcal A$; moreover, $b$
 is $\ZF$-definable from $(X,d)$ and $\mathcal A$. 
\Rp

\section{The convex topology on a normed space}
{\em In this paper,  all vector spaces that we consider are defined over the field $\IR$ of {\em real} numbers.} 
\subsection{Banach spaces}
Given a  normed space $E$ endowed with a norm $\norm$, we denote  by $B_E$ the closed unit ball  $\{x \in E : \norme{x} \le 1\}$, and by $S_E$ the unit sphere of $E$. 
The topology on $E$ associated to the norm is called the {\em strong topology}.
A {\em Banach} space is a normed space which is (Cauchy)-complete for the metric associated to the norm ({\em i.e.} every Cauchy filter of closed sets has a non-empty intersection).

\subsection{Weak topologies on normed spaces}
\subsubsection{The continuous dual $E'$ of a normed space $E$}
We endow  the vector space $E'$  of continuous linear mappings $f: E \to \IR$ 
with  the {\em dual norm} $\norm^*$, and  we call this space 
the {\em continuous dual} of the normed space $E$. We also denote by
$can: E \to E''$ the {\em canonical mapping} associating to every $x \in E$ the ``evaluating mapping'' $\tilde x : E' \to \IR$, satisfying for every $f \in E'$ the equality 
$\tilde x (f)=f(x)$. 
\subsubsection{The weak topology  $\sigma(E,E')$ on $E$}
It is the weakest topology 
$\mathcal T$ on $E$ such that elements $f \in E'$ are $\mathcal T$-continuous. 
The vector space $E$ endowed the  weak topology  is a 
 locally convex  topological vector space.   
 \subsubsection{The weak* topology  $\sigma(E',E)$ on $E'$} It is the weakest topology 
$\mathcal T$ on $E$ such that evaluating mappings $\tilde x : E' \to \IR$, $x \in E$ are $\mathcal T$-continuous. 
The vector space $E'$ endowed the  weak* topology  is a 
 Hausdorff locally convex  topological vector space.

\begin{remark} In a model of $\ZF$ where $\HB$ fails, there exists a non null (infinite dimensional) normed space $E$ such
that $E'=\{0\}$ (see \cite[Lemma~5]{Fo-Mo} or \cite{Lu-Va}). In such a model of $\ZF$,  the weak topology on $E$ is trivial with only two open sets. 
\end{remark}


\subsection{The convex topology on a normed space} \label{subsec:convex-top}
\subsubsection{Definition of the convex topology}
Since the weak topology on an infinite dimensional normed space $E$ may be trivial (in $\ZF$), we define an alternative topology on $E$, the {\em convex topology} 
(which we  introduced in \cite{Mo04}): it is the 
weakest topology $\mathcal T_c$ for which strongly closed convex subsets of $E$ are $\mathcal T_c$-closed. 
The lattice generated by strongly closed convex subsets of $E$ is called 
the {\em convex lattice} of $E$. Elements of this lattice 
are finite unions of strongly closed convex sets, so this lattice  is a basis of closed subsets of the convex topology. 
Thus the set $\mathcal C$ of strongly closed convex subsets of $E$ is a sub-basis of closed sets of the convex topology, which is  closed by finite intersection.

\begin{proposition} Let $E$ be a normed space. 
\begin{enumerate}[(i)]
\item \label{it:w-conv-strong} ``weak topology on $E$''  $\subseteq$ ``convex topology on $E$''  $\subseteq$ ``strong topology on $E$''.
\item \label{it:hb-w-conv} If $E$ satisfies  the continuous Hahn-Banach property, then the weak topology and the convex topology on $E$ are equal.
\end{enumerate}
\end{proposition}
\Pr \eqref{it:w-conv-strong} is trivial. \eqref{it:hb-w-conv} follows from the fact that if a normed space $E$ satisfies the
{\em CHB}  property, then it satisfies several classical geometric forms of the geometric Hahn-Banach property (see \cite{Do-Mo})
and in particular, every closed convex set is weakly closed.
\Rp 

 In  $\ZF \mathbf + \HB$, the weak topology and the convex topology on a normed space are equal. 
 
\subsubsection{Convex topology {\em vs} strong topology}
\begin{theorem*}[{\bf Lusternik-Schnirelmann}]  Let  $n \in \IN^*$, let  $N: \IR^n \to \IR$ be  a norm, and let $S$ be the unit sphere of $N$.
Let $a \in \IR^n$ such that $N(a)<1$. 
Denote by $s_a : S \to S$ the ``antipodal mapping'' associating to every $x \in S$ the point $y \in S$ such that $(xa) \cap S=\{x,y\}$, where 
$(xa)$ is the line generated by $x$ and $a$. 
If $C_1, \dots, C_n$ are $n$ closed subsets of $\IR^n$ such that $S \subseteq \cup_{1 \le i \le n} C_i$, then there exists
$i \in \{1..n\}$ and  $x \in S$ such that $\{x,s_a(x)\} \subseteq C_i$.
\end{theorem*}
\Pr The proof of this famous result is choiceless (see for example \cite{Mat} for $a=0$).
\Rp 

\begin{theorem} \label{theo:clos-ball} Let $E$ be a normed space which is not finite-dimensional. 
The closure of the unit sphere $S_E$ for the convex topology is the closed unit ball $B_E$.
In particular, the convex topology on $E$ is strictly contained in the strong topology. 
\end{theorem} 
\Pr Since $B_E$ is closed in the convex topology, the closure $C$ of $S_E$ in the convex topology is contained in $B_E$.
We now prove that every point $a \in E$ such that $\norme{a}<1$ belongs to $C$. 
Consider some finite set  $\{C_i : 1 \le i \le n\}$ of closed convex subsets of $E$ such that  $F:= \cup_{1\le i \le n} C_i$ contains  $S$.  We have to show that  $a \in F$. 
Let  $V$ be a vector subspace of  $E$, containing   $a$,  with dimension $\ge n$.  The  Lusternik-Schnirelmann theorem implies that for some $i_0 \in \{1..n\}$, 
 $C_{i_0} \cap V$ contains two $a$-antipodal points of   $S_E$; by convexity of $C_{i_0}$, $a \in C_{i_0}$. 
\Rp

\begin{question} The convex topology $\mathcal T_c$ on a normed space $E$ is $T_1$ ({\em i.e.} every singleton is closed). Is it  Hausdorff? Is the space
$(E,\mathcal T_c)$ a topological vector space? Is it locally convex?
Is $\mathcal T_c$ the topology associated to some family of pseudo-metrics on $E$?  Is  $\mathcal T_c$ uniformizable?
\end{question}

\section{Weak compactness in a uniformly convex Banach space}
\subsection{Uniform convexity}
\subsubsection{Strict convexity} A normed space $(E,\norm)$ is {\em strictly convex} if every segment contained  in
the unit sphere is a singleton: for every  $x, y \in S_E$, $x \neq y \Rightarrow \norme{\frac{x+y}{2}} < 1$.

\subsubsection{Uniform  convexity} This is a strong quantitative version of the strict convexity. A 
 normed space $(E,\norm)$ is {\em uniformly  convex} if for every real number $\varepsilon>0$, there exists  
$\delta >0$ such that for every  $x, y \in B_E$, 
$\big(\norme{x-y} > \varepsilon \Rightarrow \norme{\frac{x+y}{2}} < 1 - \delta \big)$.

\begin{example} Every Hilbert space is uniformly convex (see \cite[p.~190-191]{Be}). Let  $p \in ]1,+\infty[$.  If  $\mathcal B$ is a boolean algebra of subsets of a set $I$, and if 
$\nu: \mathcal B \to [0,+\infty]$ is non-null and finitely additive, then the normed space  $L^p(\nu)$ is uniformly convex (see \cite[Section~4]{Do-Mo}). In particular,
for every set $I$, the normed
space  $\ell^p(I)$ (see Section~\ref{subsec:ellp}) is uniformly convex.
\end{example}

\begin{proposition} \label{prop:uconv-small-subbasis} Given a uniformly convex normed space $E$, and denoting by $d$ the metric on $E$ given by the norm,
the family of closed convex subsets of $E$ satisfies the property
of $d$-smallness in thin $d$-crowns centered at $0_E$. 
\end{proposition}
\Pr This follows directly from the definition of uniform convexity.
\Rp

\subsection{Various weak forms of the Alaoglu theorem} \label{subsec:weak-alaoglu}
Consider the  following statements (the first two  were introduced in \cite{De-Mo} and \cite{Mo04} 
and are consequences of $\BPI$ -or rather the Alaoglu theorem-):  
\begin{itemize}
\item $ \AUc$: The closed unit ball (and thus every bounded subset which is closed in the convex topology) of a uniformly convex Banach space is compact in the convex topology.
\item $ \AH$: (Hilbert) The closed unit ball (and thus every bounded weakly closed subset) of a Hilbert space is weakly compact. 
\item $ \AHb$: (Hilbert with hilbertian basis) For every set $I$, the closed unit ball of $\ell^2(I)$ is 
weakly compact. 
\item $ \AHbf$:  For every sequence $(F_n)_{n \in \IN} $ of finite sets, the closed unit ball of 
$\ell^2(\cup_{n \in \IN} F_n)$ is weakly compact. 
\end{itemize}
Of course, $ \AUc \Rightarrow  \AH \Rightarrow  \AHb \Rightarrow  \AHbf$.  

\begin{remark} \label{rem:ah-wo} If a Hilbert space $H$ has a well orderable dense subset, then $H$ has a well orderable hilbertian basis, thus $H$ is isometrically isomorphic with some
$\ell^2(I)$ where $I$ is well orderable, and in this case, the closed ball $B_H$, which is homeomorphic with a closed subset of $[-1,1]^I$,
 is weakly compact (use Remark~\ref{rem:tycho-wo}). 
In particular, given an ordinal $\alpha$, the closed unit ball of $\ell^2(\alpha)$ (for example the closed unit ball of $\ell^2(\IN)$)
  is weakly compact.
\end{remark}

\begin{theorem} \label{theo:acd2wc-uc} 
\begin{enumerate}[(i)] 
\item \label{it:ACD2A1} $\ACD \Rightarrow  \AUc$. 
\item \label{it:A1not2ACD}  $ \AUc \not \Rightarrow \ACD$.
\item \label{it:ACD2A2} $ \AHbf \Leftrightarrow \ACDF$.    
\end{enumerate}
\end{theorem}
\Pr \eqref{it:ACD2A1}  Let $E$ be a uniformly convex Banach space. Denote by $\mathcal C$ the class of (strongly) closed convex subsets of $E$. 
Denoting  by $\mathcal T$ the convex topology on $E$, the class  $\mathcal C$  is 
a sub-basis of closed subsets of the topological space $(E,\mathcal T)$. 
Now, consider the distance $d$ associated to the norm on $E$: 
using  Proposition~\ref{prop:uconv-small-subbasis},  $\mathcal C$ satisfies the property
of $d$-smallness in thin $d$-crowns centered at $0_E$. Moreover,  the metric  space $(E,d)$ is complete, 
 closed $d$-balls belong to $\mathcal C$ and $\mathcal T$ is included in the topology associated to $d$.   Applying Theorem~\ref{theo:small-subbasis2compact}-\eqref{it:main-ACD}, it follows from $\ACD$ that the unit closed ball of $E$ (and also every bounded subset of $E$ which is closed in the convex topology) 
is compact in the convex topology of $E$. \\
\eqref{it:A1not2ACD} $\AUc \not \Rightarrow \ACD$ because $\BPI \Rightarrow \AUc$ and $\BPI \not \Rightarrow \ACD$. \\ 
\eqref{it:ACD2A2} The idea of the implication $\AHbf \Rightarrow \ACDF$ is in \cite[th.~9 p.~16]{Fo-Mo}: we sketch it for sake of completeness. Let $(F_n)_{n \in \IN}$ be a disjoint sequence of non-empty finite sets. Let us show that $\prod_{n \in \IN}F_n$ is non-empty. 
Let $I:=\cup_{n \in \IN} F_n$. Then the Hilbert spaces $H:=\ell^2(I)$ and $\oplus_{\ell^2(\IN)} \ell^2(F_n)$ are isometrically isomorph. Let $(\varepsilon_n)_{n \in \IN}$ be a sequence of $]0,1[$ such that 
$\sum_{n \in \IN} \varepsilon_n^2=1$. 
For every $n \in \IN$, let $\tilde F_n :=\{\varepsilon_n 1_{\{x\}} : x \in F_n\}$ where 
 for each $x \in F_n$, $1_{\{x\}} : F_n \to \{0,1\}$ is the  indicator of $\{x\}$; 
let $Z_n:=\{x=(x_k)_{k \in \IN} \in H : \; x_n{\restriction F_n} \in \tilde F_n\}$. Each $Z_n$  is a weakly closed  subset of the ball $B_H$ ($Z_n$ is a finite union of closed convex sets). Moreover,
the sequence $(Z_n)_{n \in \IN}$ is centered. The weak compactness  of $B_H$ implies that 
$Z:=\cap_{n \in \IN} Z_n$ is non-empty. An element of $Z$ defines an element of 
$\prod_{n \in \IN} \tilde F_n$, and thus an element of 
$\prod_{n \in \IN} F_n$ (because each $\varepsilon_n$ is $>0$). 
For the converse statement, if 
$(F_n)_{n \in \IN}$ is a sequence of finite sets, then the set $I:=\cup_{n \in \IN}F_n$ is finite or countable and
in both cases, the closed unit ball of the Hilbert space  $\ell^2(I)$ is weakly compact (see Remark~\ref{rem:ah-wo}). 
\Rp

\begin{remark} Theorem~\ref{theo:acd2wc-uc}-\eqref{it:ACD2A1} enhances our previous result 
$\DC \Rightarrow  \AUc$ which we proved in  
\cite{De-Mo},  where we left open the two questions: Does  $\ACD$ imply $\AUc$? Does $\ACD$ imply $\AH$?  
A proof of  
$\ACD \Rightarrow  \AH$ has been  found by Fremlin (see \cite[chap.~56, Section~566P]{Fr06}).
\end{remark}

\begin{figure}[h]
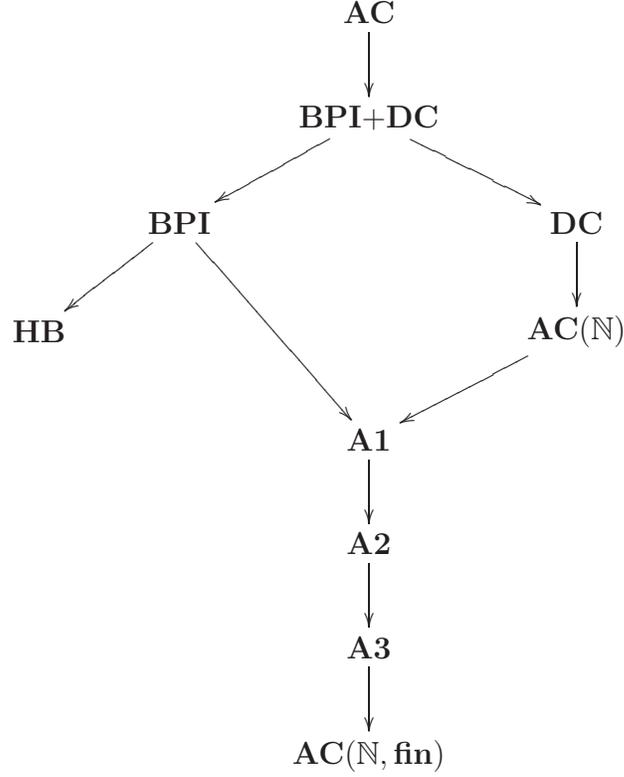

\input figure.tex
\caption{Some weak forms of $\AC$}
\label{fig:var-refl}
\end{figure}

\begin{question}
Does   $ \AH $ imply $\AUc$? Does $ \AHb$ imply $ \AH$? Does  $\ACDF$ imply $\AHb$?
\end{question}

\subsection{Convex-compactness in $\ZF$} \label{subsec:cc-dual-uc}
Given a vector  space $E$, endowed with a topology $\mathcal T$,  say that a subset $A$ of $E$ is {\em convex-compact} if, denoting by   
$\mathcal C$ the set of $\mathcal T$-closed convex subsets of $E$, $A$ is $\mathcal C$-compact; moreover, if $A$ is closely $\mathcal C$-compact, say that $A$ is 
{\em closely convex-compact}.

\begin{theorem} \label{theo:uc-ccomp-ZF} The closed unit ball of a uniformly convex Banach space is closely convex-compact in the convex topology.
\end{theorem} 
\Pr The proof is analog to the proof of Theorem~\ref{theo:acd2wc-uc}-\eqref{it:ACD2A1}, applying  Theorem~\ref{theo:small-subbasis2compact}-\eqref{it:main-ZF} instead of Theorem~\ref{theo:small-subbasis2compact}-\eqref{it:main-ACD}.
\Rp

\section{$\ACD$ and compactness  in  $[0,1]^I$} \label{sec:acd-comp}
Given a set $I$,  we endow the vector space $\IR^I$ with the  the product topology, which we denote
by  $\mathcal T_I$.

\subsection{Spaces $\ell^p(I)$, $1 \le p \le \infty$ or $p=0$} \label{subsec:ellp}
Denote by $\ell^{\infty}(I)$ the following vector space endowed with the ``sup'' norm $N_{\infty}$:
$$\ell^{\infty}(I) := \{x=(x_i)_{i \in I} : \; \sup_{i \in I} |x_i| < +\infty \}$$ 
Denote by $\ell^0(I)$ the following closed vector space of $\ell^{\infty}(I)$ endowed with the  norm $N_{\infty}$: 
$$\ell^0(I):=\{x=(x_i)_{i \in I} \in \ell^{\infty}(I)  : \; \forall \varepsilon >0 \; \exists F \in \mathcal P_f(I) \; \forall i \in I \backslash F \; |x_i| \le \varepsilon \}$$ 
For every $p \in [1,+\infty[$, denote by $\ell^p(I)$ the following vector space endowed with the $N_p$-norm:
$$\ell^p(I):=\{x=(x_i)_{i \in I} \in \IR^I : \sum_{i \in I} |x_i|^p < + \infty\}$$ 

Recall that the continuous dual of $\ell^0(I)$ is (canonically isometrically isomorphic with) $\ell^1(I)$.
Given some $p \in ]1,+\infty[$, the continuous dual of $\ell^p(I)$ is (canonically  isometrically isomorphic with) $\ell^q(I)$ where 
$q$ is the conjuguate of $p$. The following Lemma is easy:

\begin {lemma} \label{lem:weak=prod}
 Let $I$ be a set. 
\begin{enumerate}[(i)]
 \item The topology induced by $\mathcal T_I$  on the subset 
$\ell^1(I)$ ({\em resp.} $\ell^{\infty}(I)$) is included in the weak* topology $\sigma(\ell^1(I), \ell^0(I))$ of $\ell^1(I)$ 
 ({\em resp.}  the weak* topology $\sigma(\ell^{\infty}(I), \ell^1(I))$ of $\ell^{\infty}(I)$). 
 Moreover, the two topologies induce the same topology  on bounded subsets of 
 $\ell^1(I)$ ({\em resp.} $\ell^{\infty}(I)$).
\item \label{it:weak-ell0=prod} The topology induced by $\mathcal T_I$  on the subset 
$\ell^0(I)$  is included in the weak topology $\sigma(\ell^0(I), \ell^1(I))$. Moreover, the two topologies induce the same topology
 on bounded subsets of 
 $\ell^0(I)$.
\item \label{it:weak-ellp=prod} Given some $p \in ]1,+\infty[$, the topology induced by $\mathcal T_I$  on the subset 
$\ell^p(I)$ is included in the weak topology $\sigma(\ell^p(I), \ell^q(I))$ where $q$ is the conjuguate of $p$.
Moreover, the two topologies induce the same topology  on bounded subsets of  $\ell^p(I)$.
\end{enumerate}
\end {lemma}

\subsection{Closed subspaces of $[0,1]^I$ included in $\ell^p(I)$, $1 \le p <+\infty$}
Given a set $I$, denote by $[0,1]_{\sigma}^I$ the set of elements $x=(x_i)_{i \in I} \in [0,1]^I$
such that the support $\{i \in I : x_i \neq 0\}$ of $x$ is countable.  Say that a closed subset $F$ of $[0,1]^I$ is {\em Corson} if $F \subseteq [0,1]_{\sigma}^I$. 
\begin{proposition} \label{prop:sequ-comp-acd} Let $I$ be a set. Let $F$ be a closed subset of $[0,1]^I$. 
\begin{enumerate}[(i)]
\item \label{it:acdf-ell0-to-corson} If  $F \subseteq \ell^0(I)$, then $\ACDF$ implies that $F$  is Corson.
\item \label{it:acd-corson} If $F$ is Corson, then  $\ACD$ implies that $F$  is sequentially compact.
\item \label{it:seq-comp-ell0} If $F \subseteq \ell^0(I)$, then $\ACD$ implies that $F$ is sequentially compact.
\end{enumerate}
\end{proposition}
\Pr \eqref{it:acdf-ell0-to-corson} assume that $F \subseteq \ell^0(I)$. Given some  $x=(x_i)_{i \in I} \in F$, 
the support $J:=\{i \in I : \; x_i \neq 0\}$ of $x$ is a countable union of finite sets. Using $\ACDF$, $J$ is countable. Thus $F$ is Corson.\\
 \eqref{it:acd-corson} Let $(x^n)_{n \in \IN}$ be a sequence of $F$. 
For every $n \in \IN$, denote by $J_n$ the {\em support} of $x^n$: $J_n:=\{i \in I : x^n_i \neq 0\}$.
Using $\ACD$, the set $J:=\cup_{n \in \IN} J_n$ is countable. 
So $K:=[0,1]^J \times \{0\}^{I \backslash J}$ is compact and metrizable 
so $K$ is sequentially compact: extract from $(x^n)_{n \in \IN}$ a convergent subsequence $(x^n)_{n \in A}$ where $A$ is some infinite subset of $\IN$. \\
\eqref{it:seq-comp-ell0}  Use \eqref{it:acdf-ell0-to-corson} and \eqref{it:acd-corson}. 
\Rp

\begin{corollary}
Let $F$ be a  closed subset of $[0,1]^I$. If there exists $p \in [1,+\infty[$ such that 
$F$ is a bounded subset of $\ell^p(I)$, then $\ACD$ implies that $F$  is compact in $[0,1]^I$.
\end{corollary}
\Pr Let $r \in ]p,+\infty[$. Then   $N_r \le N_p$, so $F$ is a bounded subset of $\ell^r(I)$.  
Since $F$ is weakly closed and bounded in $\ell^r(I)$, 
 Theorem~\ref{theo:acd2wc-uc}-\eqref{it:ACD2A1} implies that, using $\ACD$,   $F$ is 
 compact in the weak topology $\sigma(\ell^r(I),\ell^{r'}(I))$ where $r'$ is the conjuguate of $r$. It follows from 
 Lemma~\ref{lem:weak=prod}-\eqref{it:weak-ellp=prod} that $F$ is compact
 for the topology $\mathcal T_I$.  
\Rp

\begin{question}
 What is the power of the statement {\em ``The closed unit ball of $\ell^2(\IR)$ is weakly compact''}?
This statement is a consequence of $\ACD$. Are there models of $\ZF$ which do not satisfy this statement? 
\end{question}

\section{$\DC$ and compactness in  $[0,1]^I$} \label{sec:dc-and-comp}
\subsection{Eberlein's criterion of compactness}
\subsubsection{$\vartheta$-sequences}
Let $E$ be a normed space, and denote by $d$ the metric given by the norm on $E$. 
Given a subset $F$ of $E$, and some real number $\vartheta>0$, a  {\em $\vartheta$-sequence of $F$} is a sequence  $(a_n)_{n \in \IN}$ of  $F$ 
satisfying for every $n \in \IN$: 
$$d\big({\spanv}\{a_i : i <n\}, \conv \{a_i : i \ge n\} \big) \ge \vartheta$$
In \cite{Mo04}, we proved in $(\ZF \mathbf + \DC)$ the following result:
\begin{theorem*}[$\DC$] Let  $E$ be a Banach space. Denote by  $\mathcal T$ the convex topology on $E$.  
Let   $F$ be a   {\em convex} subset of $E$ which is  $d$-bounded and $\mathcal T$-closed. If $F$ is not $\mathcal T$-compact,
then there exists some real number  $\vartheta >0$ and a  $\vartheta$-sequence of $F$. 
\end{theorem*}
If we delete the hypothesis ``$F$ is convex'', our next result allows us to build  in $(\ZF \mathbf + \DC)$ ``pseudo $\vartheta$-sequences''.
Say that a sequence $(a_n)_{n \in \IN}$ of $F$ is a 
 {\em pseudo $\vartheta$-sequence} if for every $n \in \IN$:
$$d\big({\spanv}\{a_i : i <n\}, (F \cap \conv \{a_i : i \ge n\}) \big) \ge \vartheta$$

\subsubsection{Building pseudo-sequences with $\DC$}
We first recall the following result for saturating filters {\em w.r.t.} some numeric constraint:
\begin{proposition*}[$\DC$]
Let $E$ be a set, let  $\mathcal L$ be a lattice of subsets of $E$, with smallest element  
$\varnothing$ and greatest element  $E$. Let  $\rho : \mathcal L \to \IR_+$ be some mapping. 
let  $\tilde \rho : \mathcal P(\mathcal L) \to \IR_+$ be the mapping associating to every subset 
$\mathcal A$ of  $\mathcal L$ the real number  $\inf \{\rho(A) : A \in \mathcal A\}$.
Let $\mathcal F$ be a filter of  $\mathcal L$. 
Then there exists a filter  $\mathcal G$ of  $\mathcal L$ including  $\mathcal F$
such that  
$$\tilde \rho(\mathcal G) = \tilde \rho(\mathcal S_{\mathcal L}(\mathcal G))$$ 
\end{proposition*}
\Pr See   \cite{Mo04}. 
\Rp 

\begin{remark} The previous Proposition is easy to prove in $\ZFC$: consider a maximal filter of $\mathcal L$
including $\mathcal F$.  
\end{remark}

\begin{notation} Given a metric space $(X,d)$ and    some subset $A$ of $X$, we denote by $d_A$ the ($1$-Lipschitzian hence)
continuous mapping $d_A: X \to \IR$ associating to every $x \in X$ the real number $d(x,A):=\inf\{d(x,a): a \in A\}$.
Moreover, given  some real number $\vartheta >0$,
we denote by $A_{\vartheta}$ the following closed subset of $X$:
$$A_{\vartheta} := \{x \in X : \; d(x,A) \le \vartheta\}$$
\end{notation}

Notice that if $A$ is a convex subset of a normed space, then for every $\vartheta >0$ the set $A_{\vartheta}$ is convex (because
the mapping $d_A$  is convex).

\begin{theorem}[$\DC$] \label{theo:jseq-ferme}
Let $E$ be a Banach space. Let $d$ be the distance given by the norm on $E$. Let  $\mathcal T$
be the convex topology on $E$. Let  $F$ be a $d$-bounded subset of  $E$,  which is  $\mathcal T$-closed (thus $d$-closed, thus  $d$-complete).
 If  $F$ is not   $\mathcal T$-compact, then there exists some real number $\vartheta >0$, and a sequence  
$(a_n)_{n \in \IN}$ of $F$, such that for every $n \in \IN$: 
$$d\big(\spanv\{a_k : k \le n\}, (\overline{\conv}^{\mathcal T}\{a_k : k>n\} \cap F)\big) \ge \vartheta$$
\end{theorem}
\Pr Let   $\mathcal C$ be the set of $\mathcal T$-closed ({\em i.e.} $d$-closed) convex subsets of $E$. 
Let  $\mathcal L_c$ be the lattice generated by   $\mathcal C$. 
Let  $\mathcal L_1$ be the lattice induced by  $\mathcal L_c$ on   $F$: 
$\mathcal L_1 =  \{A \cap F : \; A  \in \mathcal L_c\}$.
Since $F$ is not  $\mathcal T$-compact, let  $\mathcal F$ be a filter  of $\mathcal L_1$ containing $F$ such that $\cap \mathcal F = \varnothing$.
Let  $\rho$ be the  ``diameter'' function ({\em w.r.t.} the distance $d$), which is defined for   
$d$-bounded subsets of  $E$, and in particular on  $\mathcal L_1$. 
Using the previous Proposition for the ``diameter'' function $\rho$,  $\DC$ implies the existence of a filter  
  $\mathcal G$  of $\mathcal L_1$ including  $\mathcal F$ such that  
$$r:=\tilde \rho(\mathcal S_{\mathcal L_1}(\mathcal G))=\tilde \rho(\mathcal G)$$
If  $r=0$, then, since  the metric space $(F,d)$ is complete, $\cap \mathcal G$ is a singleton  $\{a\}$,
and $a \in \cap \mathcal F$: contradictory! Thus $r>0$. Let  $0<\vartheta<r$.   
We will now build a sequence  $(K_n)_{n \in \IN}$ of   $\mathcal C$, and a sequence  $(a_n)_{n \in \IN}$ of $F$,
such that for every  $n \in \IN$, 
\begin{equation} \label{equ:pseudo-james} 
a_n \in \cap_{i \le n} K_i \text{ and }  (\spanv \{a_i : i < n\})_{\vartheta} \cap (K_n \cap F) = \varnothing
\end{equation}
 It will follow that  for every $n \in \IN$,
 $$d\big(\spanv\{a_k : k < n\}, (\overline{\conv}^{\mathcal T}\{a_k : k \ge n\} \cap F)\big) \ge \vartheta$$
\par $\bullet$ $B(0,\vartheta)  \notin \mathcal S_{\mathcal L_1}(\mathcal G)$: thus there exists  
 $G \in \mathcal G$ satisfying   $B(0,\vartheta) \cap G = \varnothing$; since  $G \in \mathcal L_1$, $G$ is of the form  
$\cup_{i \in I}  (C_i \cap F)$ where $I$ is finite and each  $C_i$ belongs to $\mathcal C$; using Section~\ref{subsubsec:stat}-\eqref{it:stat2},
let    $i \in I$ such that  
 $(C_i \cap F) \in \mathcal S(\mathcal G)$.
let  $K_0$ be the convex set  $C_i$.  Let  $\mathcal G_0$ be the filter of $\mathcal L_1$ generated by  
 $\mathcal G$ and $K_0$. Let    $a_0 \in K_0 \cap F$.
\par  $\bullet$  $(\IR a_0)_{\vartheta} \notin \mathcal S_{\mathcal L_1}(\mathcal G_0)$: let   $G \in \mathcal G_0$ such that  
$(\spanv \{a_0\})_{\vartheta}  \cap G = \varnothing$.   
since   $G \in \mathcal L_1$,  $G$ is of the form  
$\cup_{i \in I} (C_i \cap F)$ where $I$ is finite and each $C_i$ belongs to $\mathcal C$; let    $i \in I$ such that  
$(C_i \cap F) \in \mathcal S(\mathcal G_0)$.
Let $K_1$ be the set $C_i$.  Let $\mathcal G_1$ be the filter of  $\mathcal  L_1$ generated
by  $\mathcal G_0$ and  $K_1$. Let  $a_1 \in  (K_0 \cap K_{1} \cap F)$. 
\par  $\bullet$  $(\spanv \{a_0,a_1\})_{\vartheta} \notin \mathcal S_{\mathcal L_1}(\mathcal G_1)$: 
 let  $K_2 \in \mathcal C$ such that  $K_2 \in \mathcal S(\mathcal G_1)$  and  
$(\spanv \{a_0,a_1\})_{\vartheta} \cap (K_2 \cap F) = \varnothing$.  
Let  $\mathcal G_2$ be the filter of  $\mathcal  L_1$ generated by  $\mathcal G_1$ and $K_2$. 
Let  $a_{2} \in (K_{0} \cap K_{1} \cap K_2) \cap F$. 
\par  $\bullet$  \dots
\par Using $\DC$, we construct  a sequence  $(a_n,C_n)_{n \in \IN}$ of $F \times \mathcal C$ satisfying \eqref{equ:pseudo-james}.
\Rp

\begin{corollary} \label{cor:comp-top-c}
Let $E$ be a Banach space. Let $d$ be the metric given by the norm on  $E$. 
Let $F$ be a $d$-bounded subset of $E$ which is closed for the convex topology $\mathcal T$ on $E$.
Consider the three following statements:
\begin{enumerate}[(i)]
\item \label{it:comp1} $F$ is compact for the convex topology $\mathcal T$.
\item \label{it:comp2} $F$ is sequentially compact for $\mathcal T$.
\item \label{it:comp3} For every $\vartheta >0$, $F$ does not contain any pseudo $\vartheta$-sequence.
\end{enumerate} 
Then \eqref{it:comp1} $\Rightarrow$ \eqref{it:comp2} $\Rightarrow$ \eqref{it:comp3}. Moreover, in $(\ZF \mathbf + \DC)$,
\eqref{it:comp3} $\Rightarrow$ \eqref{it:comp1}.
\end{corollary}
\Pr \eqref{it:comp1} $\Rightarrow$ \eqref{it:comp2} It is sufficient to prove this implication when $E$ is a  {\em separable} Banach space.
In this case, the normed space satisfies the {\em CHB} property (see Section~\ref{subsubsec:HB}), thus the convex topology 
 $\mathcal T$ and the weak topology on $E$ are
equal. Moreover, there is a norm $N$ on $E$ which induces on the closed unit ball of $E$ a topology which is included in the weak topology of $E$ (see
for example \cite[Lemme~I.4 p.~2]{Li-Qu}).
This implies that the topology given by the norm $N$ and the topology  $\mathcal T$ are equal on $K$; thus the 
 $\mathcal T$-compact space $K$ is metrisable whence $K$ is sequentially compact.  \\
\eqref{it:comp2} $\Rightarrow$ \eqref{it:comp3}  Assume that the subset $K$ is sequentially compact in the topology $\mathcal T$. 
Let $\vartheta >0$. Seeking for a contradiction, assume that $F$ has a pseudo $\vartheta$-sequence $(a_n)_{n \in \IN}$.
Extract some sequence  $(a_n)_{n \in A}$ which converges to some $l \in F$ in the  topology $\mathcal T$. 
Then, for every  $n \in \IN$,
$l \in \overline{\conv}^{\mathcal T} \{a_i : i \ge n\}$. 
Let  $V:=\spanv \big(\{a_i : i \in \IN\} \cup \{l\}\big)$. 
 Let  $(u_n)_n$ be a convex block-sequence of  $(a_n)_{n \in A}$
which strongly converges to  $l$. Since the normed space $V$ is separable, for each $n \in \IN$, choose some $f_n$ in the unit
sphere of $V'$  such that $f_n$ is null on  $\{a_i : i < n\}$ and  $f_n(l) \ge \vartheta$.  Let 
$n_0 \in \IN$ such that  
$d(u_{n_0},l) < \frac{\vartheta}{2}$. Let  $N \ge n_0$ such that  $u_{n_0} \in
\spanv \{a_i : i < N\}$; then 
$f_N(u_{n_0}) =0$ and   $f_N(l) \ge \vartheta$ thus,  since $\norme{f_N}=1$, 
$d(l, u_{N_0}) \ge \vartheta$: this is contradictory! \\
The implication  \eqref{it:comp3} $\Rightarrow$ \eqref{it:comp1} holds in $\ZF \mathbf + \DC$ thanks to  Theorem~\ref{theo:jseq-ferme}.
\Rp

\subsection{Closed subsets of $[0,1]^I$ included in $\ell^0(I)$}
\begin{corollary} \label{cor:comp-dc-ell0}
Let $F$ be a  closed subset of $[0,1]^I$. Assume that  
$F$ is a (bounded) subset of $\ell^0(I)$.  Then  $\DC$ implies that $F$  is compact.
\end{corollary}
\Pr The normed space  $\ell^0(I)$ satisfies the {\em CHB} property (see Section~\ref{subsubsec:HB}), thus 
the weak topology $\sigma(\ell^0(I),\ell^1(I))$ and the convex topology $\mathcal T$ on $\ell^0(I)$ are equal on $\ell^0(I)$.
Since $F$ is a bounded subset of $\ell^0(I)$, the topology $\mathcal T$ and the product topology $\mathcal T_I$ induce the same topology  on $F$ (see Lemma~\ref{lem:weak=prod}-\eqref{it:weak-ell0=prod}).
The subset $F$ of $\ell^0(I)$ is bounded and $\mathcal T$-closed; using Proposition~\ref{prop:sequ-comp-acd}-\eqref{it:seq-comp-ell0},
$F$  is sequentially compact for the topology $\mathcal T_I$ {\em i.e.}  for  $\mathcal T$.
Using $\DC$, Corollary~\ref{cor:comp-top-c} implies that  $F$ is compact for $\mathcal T$ {\em i.e.} for $\mathcal T_I$.
\Rp

\begin{question} Let $I$ be a set and $F$ be some closed subset of $[0,1]^I$. 
If $F \subseteq \ell^0(I)$, does $\ACD$ imply that $F$ is compact? 
\end{question}

\begin{question}
 Let $F$ be a closed subset of $[0,1]^I$ which is contained in $\IR^{(I)}$ (the vector subspace of elements 
$x \in \IR^I$ which have a {\em finite} support). Then $F \subseteq \ell^0(I)$ thus, using Corollary~\ref{cor:comp-dc-ell0},
$\DC$  implies that  $F$ is compact. Does $\ACD$ imply that $F$ is compact? 
\end{question}

\begin{question} Let $I$ be a set and $F$ be some closed subset of $[0,1]^I$. 
\begin{enumerate}[(i)] 
\item If $F$  is Corson, 
does $\DC$ imply that $F$ is compact? Does $\ACD$ imply that $F$ is compact? ($\DC_{\aleph_1}$ implies that
$F$ is compact).
 \item More generally, which closed subsets of $[0,1]^I$ can be proved compact in  $\ZF \mathbf + \DC$?
in  $\ZF \mathbf + \ACD$? 
\end{enumerate}
\end{question}

\bibliographystyle{abbrv}
\bibliography{../biblio}
\end{document}